\documentclass[11pt,a4paper]{article}
\usepackage[titletoc]{appendix}
\usepackage[labelsep=period]{caption}
\usepackage{bm}
\usepackage{fullpage}
\usepackage{epsfig}
\usepackage{graphics}
\usepackage{latexsym}
\usepackage{amsmath}
\usepackage{amsfonts}
\usepackage{amssymb}
\usepackage{mathrsfs}
\usepackage{pifont}
\usepackage{yhmath}
\usepackage{fullpage,graphicx}
\usepackage{enumitem}
\usepackage{subcaption}
\usepackage{float}
\usepackage{mathtools}
\usepackage{xcolor}
\usepackage{bbm}
\usepackage{dsfont}
\usepackage{eufrak}
\usepackage{wasysym}
\usepackage{underscore}
\usepackage{epstopdf}
\usepackage[T1]{fontenc}
\usepackage{amsthm}
\usepackage{hyperref}
\usepackage{microtype}
\usepackage{lmodern}
\usepackage[margin=25mm]{geometry}
\newtheorem{theorem}{Theorem}[section]

\newtheorem{lemma}[theorem]{Lemma}

\newcommand{\G}{\mathcal G}

\newcommand{\R}{\mathbb R}
\newcommand{\Z}{\mathbb Z}

\DeclareMathAlphabet{\mathpzc}{OT1}{pzc}{m}{it}
\title{\bf List colouring of circular-arc graphs}
\author{\vspace{2mm}  Rongchuan Tao\thanks{Corresponding author. E-mail: rctao@hku.hk.}
\\
Department of Mathematics\\University of Hong Kong \\ Hong Kong, China\smallskip\\
}
\begin{document}
\date{}
\maketitle

\begin{abstract}
A circular-arc graph takes a set of nonempty closed arcs on a circle as its vertex set, and two vertices are adjacent when their corresponding arcs intersect on the circle. We show the chromatic number of any circular-arc graph equals its list chromatic number. Our methodology is an application of Schauz's coefficient formula within the Alon-Tarsi and Combinatorial Nullstellensatz framework.

\vskip 4mm


\noindent {\bf OR/MS subject classification.} Primary: Programming/graphs.


\end{abstract}

\section{Introduction}

A \emph{circular-arc graph} takes a set of nonempty closed arcs on a circle as its vertex set, and two vertices are adjacent when their corresponding arcs intersects on the circle. The class of circular-arc graphs serves as a natural generalization of the class of interval graphs and arises in a variety of applications. For example, the wavelength-assignment problem asks the minimum number of distinct wavelengths assigned to each communication path on a circle so that paths sharing a fiber link use different wavelengths, which is exactly the colouring problem of circular-arc graphs~\cite{Valencia2003}. While the colouring problem for interval graphs can be solved efficiently by greedy algorithm, it is NP-hard for circular-arc graphs in general as shown by Garey, Johnson, Miller and Papadimitriou~\cite{Garey1980}. Indeed, circular-arc graphs can contain induced odd cycles, which implies such graphs need not be perfect.

The investigation of the colouring problem for circular-arc graphs has a long history. Tucker~\cite{Tucker1975} gave a simple approximation algorithm of approximation ratio 2, which was improved to 5/3 by Shih and Hsu~\cite{ShihHsu1990}. Karapetian~\cite{Karapetian1980} proved Tucker's conjectured bound that the chromatic number is at most 3/2 times the size of maximum clique. Valencia-Pabon~\cite{Valencia2003} subsequently refined the analysis of Tucker's algorithm using the maximum number of arcs through a point and the minimum number needed to cover the circle. Orlin, Bonuccelli and Bovet~\cite{Orlin1981} gave a polynomial algorithm for proper circular-arc graphs. Chen, Hu and Zang~\cite{Chen2005} gave an $O(n^2m)$ combinatorial algorithm for perfect circular-arc graphs.

Despite the abundance in the colouring results for circular-arc graphs, there are far less results concerning the list chromatic number of circular-arc graphs. In this paper we prove that for any circular-arc graph, its chromatic number is equal to its list chromatic number. Therefore, the list chromatic number of any circular-arc graph is at most 3/2 times the size of maximum clique. Our proof first reduces the problem to ``maximal'' circular-arc graphs. Then we construct a class $\G$ of graphs that contains all such maximal graphs. Finally, we use polynomial method to show every graph in $\G$ satisfies the property. This polynomial approach can be traced back to Alon and Tarsi~\cite{AlonTarsi1992}, who related coefficients of graph polynomials to signed counts of Eulerian subdigraphs. Alon's Combinatorial Nullstellensatz~\cite{Alon1999} provides a general
criterion for obtaining nonzero polynomial evaluations from a nonzero coefficient of certain monomial. The coefficient formula of Schauz~\cite[Theorem~3.3]{Schauz2008} expresses such a coefficient as a weighted sum of evaluations on cartesian product of finite sets.

All graphs below are finite and simple. Let $\G$ consists of the graphs $G$ such that any $v\in V(G)$ is contained in a clique $Q_v$ of size $\chi(G)$ and for any $uv\in E(G)$, either $u\in Q_v$ or $v\in Q_u$.

For any integer $n$, let $[n]:=\{1,2,\ldots,n\}$. Any arc $X$ in a circle $S$ is defined by its two end-points, listed in the clockwise order as $X^{\vdash}$ and $X^\dashv$. A \emph{k-colouring} of a graph $G$ is a map $c:V \to [k]$ such that $c(u)\neq c(v)$ for any adjacent vertices $u,v\in V(G)$. The chromatic number $\chi(G)$ is the least number $k$ such that $G$ admits a $k$-colouring. A \emph{list assignment} assigns each vertex $v$ a finite set $L(v)$ of positive integers. An \emph{$L$-colouring} is a proper colouring $c$ with $c(v)\in L(v)$ for every vertex. A graph $G$ is \emph{$k$-choosable} if for any list assignment $L$ with $|L(v)|=k$ for any $v\in V(G)$, the graph $G$ admits an $L$-colouring. The list chromatic number $\chi_L(G)$ is the least number $k$ such that $G$ is $k$-choosable. Our main result is as follows.

\begin{theorem}\label{main}
Any circular-arc graph $G$ satisfies $\chi(G)=\chi_L(G)$.
\end{theorem}

\section{Reduction}

In this section we reduce the main theorem to the following form. 

\begin{theorem}\label{reduction}
Any graph $G\in \G$ satisfies $\chi(G)=\chi_L(G)$.
\end{theorem}

\begin{proof}[Proof of theorem 1.1 by assuming theorem 2.1]

Let $G$ be a circular arc graph with $\chi(G)=k$. Fix a proper colouring $f:V(G)\to [k]$ and an arc representation $\{X_v\}_{v\in V(G)}$ on a circle $S$ with a fixed clockwise direction.

Notice that extending an arc of $G$ in the clockwise direction does not increase the chromatic number of $G$ unless the new end exceeds an endpoint of an arc in the same colour class. Since an extended arc may intersect more other arcs in $G$, repeatedly performing such operation results in a circular arc graph with potentially more edges while preserving $k$-chromatic. This is exactly the rationale behind the construction of maximal graphs.

For distinct $p,q\in S$, let $[p,q)$ be the subset of $S$ that corresponds to the clockwise arc from $p$ to $q$, including $p$ and excluding $q$. Let $K$ be a colour class of $G$. For each vertex $v\in K$, let $p_v:=X_v^\vdash$, i.e. the clockwise start of $X_v$. Now we define $T_v=[p_v,p_{v'})$ if $|K|\geq 2$, where $v'$ is the vertex in $K$ such that $p_{v'}$ immediately follows $p_v$ in the clockwise order. Otherwise $|K|=1$ and we define $T_v=S$, where $K=\{v\}$.

Since $v,v'$ are in the same colour class, we have $p_{v'}\notin X_v$ and hence $X_v\subsetneq T_v$. By definition, the set $\{T_v:v\in K\}$ gives a partition of $S$. We have therefore constructed extended arc $T_v$ for each original arc $X_v$ with one same end $p_v$. In particular, they satisfy
\begin{equation}\label{eq:partitions}
 X_v\subseteq T_v,\qquad p_v\in T_v,
 \qquad
 \{T_v:f(v)=i\}\text{ partitions }S\quad(1\le i\le k).
\end{equation}

Now let $H$ be the intersection graph of the sets $\{T_v:v\in V(G)\}$. Then $V(H)=V(G)$ and the containments in \eqref{eq:partitions} imply $E(G)\subseteq E(H)$. Thus if $H$ is $k$-choosable, then $G$ is also $k$-choosable. Since our construction allows no extension of an arc to overpass a starting point of another arc in the same colour class, we have $\chi(H)=k$. Thus it remains to show $H\in \G$ by Theorem \ref{reduction}.

For every $v\in V(H)$, define $Q_v=\{w\in V(G):p_v\in T_w\}$. Then exactly one set of each original colour contains $p_v$, by~\eqref{eq:partitions}. Hence $|Q_v|=k$. These sets all meet at $p_v$, so $Q_v$ is a clique of $H$ containing $v$. So in order to show $H\in \G$, it remains to verify every edge $e\in E(H)$ is covered by $Q_v$ for some $v\in e$. Suppose $uv\in E(H)$, so $T_u\cap T_v\neq\varnothing$. If one of them, say $T_u$, is the entire circle $S$, then $p_v\in T_u$ and $u\in Q_v$. Thus we may assume both of $T_u$ and $T_v$ are half-open arcs. Then by definition we have $T_u=[p_u,p_u')$ and $T_v=[p_v,p_v')$ for some $u',v'\in V(G)$. Now $T_u\cap T_v\neq\varnothing$ implies either $p_u\in T_v$ or $p_v\in T_u$. Therefore either $u\in Q_v$ or $v\in Q_u$. It follows that $H\in\G$. 
\end{proof}

\section{Proof of the theorem \ref{reduction}}
We first recall the following special case of Schauz's coefficient formula~\cite[Theorem~3.3(i)]{Schauz2008}. We include a proof for completeness.

\begin{lemma}\label{lem:coefficient}
Let $V$ be a nonempty finite set, let $k\ge1$, and
$M:=\prod_{v\in V}x_v^{k-1}$. For each $v$, let $B_v\subseteq\R$ have exactly $k$ distinct elements. If a polynomial $F$ has total degree at most $|V|(k-1)$, then the coefficient of $M$ in $F$ is
\begin{equation}\label{eq:coefficient}
\sum_{c\in\prod_{v\in V}B_v}
 \frac{F((c(v))_{v\in V})}
 {\displaystyle\prod_{v\in V}
    \prod_{a\in B_v\setminus\{c(v)\}}(c(v)-a)}.
\end{equation}
If $M$ does not appear in $F$, then the coefficient of $M$ in $F$ is 0. The member $c$ of the product chooses one value $c(v)\in B_v$ for every $v$. 
\end{lemma}

\begin{proof}
Suppose $B\subseteq \R$, $n\in \Z$ with $|B|=k$ and $0\le n\le k-1$. The following polynomial identity holds since both polynomials have degree at most $k-1$ and take the same value at the $k$ points of $B$: each product of the right hand side is a base polynomial of Lagrange interpolation over $B$, which attains 1 at $b$ and 0 at $B\setminus\{b\}$.
\[
 x^n=\sum_{b\in B}b^n
       \prod_{a\in B\setminus\{b\}}\frac{x-a}{b-a}
\]
Comparing coefficients of $x^{k-1}$ for each $0\leq n\leq k-1$ gives
\begin{equation}\label{eq:univariate}
 \sum_{b\in B}\frac{b^n}
 {\displaystyle\prod_{a\in B\setminus\{b\}}(b-a)}
 =\begin{cases}
    0,&0\le n<k-1,\\
    1,&n=k-1.
   \end{cases}
\end{equation}

Suppose $F$ contains a monomial $N=\prod_v x_v^{n_v}$ with coefficient $t$. Then $\sum_v n_v\le |V|(k-1)$ since $F$ has total degree at most $|V|(k-1)$. Then the contribution of $N$ to (\ref{eq:coefficient}) reads
\[
\sum_{c\in\prod_{v\in V}B_v}
 \frac{t\cdot\prod_{v\in V}c(v)^{n_v}}
 {\displaystyle\prod_{v\in V}
    \prod_{a\in B_v\setminus\{c(v)\}}(c(v)-a)}=
t\cdot \left(\sum_{c\in\prod_{v\in V}B_v}\prod_{v\in V}
 \frac{c(v)^{n_v}}
 {\displaystyle\prod_{a\in B_v\setminus\{c(v)\}}(c(v)-a)} \right)
\]
Notice the above expression factors as
\[
t\cdot \prod_{v\in V}\left(
   \sum_{b\in B_v}
   \frac{b^{n_v}}{\displaystyle\prod_{a\in B_v\setminus\{b\}}(b-a)}
 \right)
\]
Since $\sum_v n_v\le |V|(k-1)$ and by~\eqref{eq:univariate}, the above expression is $t$ if and only if $n_v=k-1$ for any $v\in V$ and otherwise zero by~\eqref{eq:univariate}. Thus any monomial in $F$ distinct from $M$ contributes zero to (\ref{eq:coefficient}), and if $F$ contains the monomial $M$, then it contributes $t$ to (\ref{eq:coefficient}). This proves the equality.
\end{proof}

Now we are ready to prove theorem \ref{reduction}. Let $H\in \G$ with $\chi(H)=k$. We show that $H$ is $k$-choosable. Write $V=V(H)$, $n=|V|$, and define polynomial $P$ on variables $\{x_v:v\in V\}$ over $\R$ by
\begin{equation}\label{eq:polynomial}
 P(\boldsymbol{x})=
  \prod_{v\in V}\ \prod_{w\in Q_v\setminus\{v\}}(x_v-x_w),
 \qquad M=\prod_{v\in V}x_v^{k-1}.
\end{equation}
Since any clique $Q_v$ has size $k$, there are $k-1$ factors indexed by each vertex. Thus $\deg P=n(k-1)=\deg M$. Since each $Q_v$ is a clique containing $v$, every factor corresponds to an edge of $H$. Conversely by assumption, for every edge $uv$ of $H$, either $(x_u-x_v)$ or $(x_v-x_u)$ is a factor in $P$. Thus for any map $c:V\to\R$, regarded as a colouring of $H$, we have 
\begin{equation}\label{eq:proper}
  P(c)\ne0\quad\Longleftrightarrow\quad
  c\text{ is a proper colouring of }H.
\end{equation}

Apply Lemma~\ref{lem:coefficient} by taking $B_v=\{1,\ldots,k\}=:B$ for every vertex $v$. If a colouring $c:V\to B$ is proper, the $k$ vertices in each $Q_v$ receive all $k$ labels exactly once. Hence for any $v\in V$ we have
\[
 \prod_{w\in Q_v\setminus\{v\}}(c(v)-c(w))
 =\prod_{a\in B\setminus\{c(v)\}}(c(v)-a).
\]
It follows that 
\begin{equation*}
 \frac{P(c)}{\displaystyle\prod_{v\in V}\prod_{a\in B_v\setminus\{c(v)\}}(c(v)-a)}=
 \frac{\prod_{v\in V}\prod_{w\in Q_v\setminus\{v\}}(c(v)-c(w))}{\displaystyle\prod_{v\in V}\prod_{a\in B_v\setminus\{c(v)\}}(c(v)-a)}=1.
\end{equation*}

Conversely, if $c:V\to B$ is improper, then $P(c)$ is zero. By (\ref{eq:coefficient}), the coefficient of $M$ in $P$ is exactly the cardinality of the set $\{c:V\to B:\ c\text{ properly colours }H\}$. This coefficient of $M$ in $P$ is therefore positive since that set is nonempty by the $k$-colourability of $H$.

Now take any list assignment $L$ with each $|L(v)|=k$. Apply Lemma~\ref{lem:coefficient} by setting $B_v=L(v)$ for each $v$. The coefficient on the left is still the positive number in~\eqref{eq:coefficient}. All denominators are nonzero. Consequently at least one choice $c(v)\in L(v)$ has $P(c)\ne0$. By~\eqref{eq:proper}, this choice is a proper list colouring of $H$.

\bibliographystyle{plain}
\bibliography{introduction-references}
\end{document}